\documentclass[12pt]{article}
\usepackage{amsmath}
\usepackage{amstext}
\usepackage{amsfonts}
\usepackage{amsthm}
\usepackage{amscd}
\numberwithin{equation}{section}

\swapnumbers
\theoremstyle{plain}
\newtheorem{thm}{Theorem}[section]
\newtheorem{thmi}{Theorem}
\newtheorem{lem}[thm]{Lemma}
\newtheorem{prop}[thm]{Proposition}
\newtheorem{cor}[thm]{Corollary}
\theoremstyle{remark}
\newtheorem{rem}[thm]{Remark}

\def\mod0{{\mathcal SU}_C(2,{\mathcal O})}
\def\modk{{\mathcal SU}_C(2,\omega)}
\def\modalpha{{\mathcal SU}_C(2,\alpha)}
\def\modkalpha{{\mathcal SU}_C(2,\omega \alpha)}

\def\cO{{\mathcal O}}

\def\D{\Delta}
\def\cL{{\mathcal L}}

\def\map#1{\ \smash{\mathop{\longrightarrow}\limits^{#1}}\ }

\def\bP{\mathbb{P}}
\def\bC{\mathbb{C}}

\def\TT{\mathbb{T}}
\def\T{\Theta}
\def\Lalpha{\mathcal{L}_\alpha}
\def\Lomegaalpha{\mathcal{L}_{\omega \alpha}}
\def\tC{\tilde{C}}
\def\tp{\tilde{p}}
\def\tq{\tilde{q}}
\def\tE{\tilde{E}}
\def\tt{\tilde{t}}
\def\ts{\tilde{s}}

\def\tTheta{\tilde{\Theta}}
\def\tT{\tilde{\Theta}}
\def\inj{\hookrightarrow}

\def\sym{\mathrm{Sym}}
\def\ker{\mathrm{ker} \:}
\def\pic{\mathrm{Pic}}
\def\mult{\mathrm{mult}}

\def\dim{\mathrm{dim} \:}
\def\det{\mathrm{det} \:}

\def\im{\mathrm{im} \:}

\def\lra{\longrightarrow}

\def\ra{\rightarrow}

\def\lms{\longmapsto}

\def\cc{\mathbb{C}}
\def\pp{\mathbb{P}}

\def\cD{\mathcal{D}}

\begin{document}

\title{Some properties of second order theta functions on Prym varieties}
\author{E. Izadi and C. Pauly}
\maketitle

\begin{abstract}
Let $P \cup P'$ be the two component Prym variety associated to an
\'etale double cover $\tC \ra C$ of a non-hyperelliptic curve
of genus $g \geq 6$ and let $|2\Xi_0|$ and $|2\Xi_0'|$ be the linear
systems of second order theta divisors on $P$ and $P'$ respectively. 
The component $P'$ contains
canonically the Prym curve $\tC$. We show that the base locus
of the subseries of divisors containing $\tC \subset P'$ is
exactly the curve $\tC$. We also prove canonical isomorphisms
between some subseries of $|2\Xi_0|$ and $|2\Xi_0'|$ and some 
subseries of second order theta divisors on the Jacobian of $C$. 
\end{abstract}

\begin{center}
{\Large \sc Introduction}
\end{center}

Let $C$ be a curve of genus $g\geq 5$ with an \'etale double cover $\pi:
\tilde{C} \ra C$. Let $Nm : Pic (\tC)\ra Pic(C)$ be the norm map.
Consider the Prym varieties
$$ Nm^{-1}(\cO) = P \cup P'
$$
which are characterized by the facts that $\cO \in P$, $\cO\not\in
P'$. Let $\sigma :\tC\ra\tC$ be the
involution of the cover $\pi :\tC\ra C$. The curve $\tC$ admits a
natural embedding in $P'$ given by the morphism
\begin{eqnarray*}
i :\tC & \lra & P'\\
\tp & \lms & \cO_{\tC } (\tp -\sigma\tp ).
\end{eqnarray*}
A symmetric Riemann theta divisor $\tT_0$ on the Jacobian $J\tC$ of
$\tC$ induces twice a symmetric principal polarization $\Xi_0$ on $P$
(resp. $\Xi_0'$ on $P'$). Let $\Gamma_{\tC }$ be the space of sections
of $\cO_{ P'} (2\Xi_0')$ vanishing on the image of $i$. In his work on
the Schottky problem, Donagi proved in \cite{donagi} (Lemma 4.8 page
597) that the base locus Bs$(\bP\Gamma_{\tC })$ of $\pp \Gamma_{\tC }$
is $i(\tC)$ for a Wirtinger cover $\pi :\tC\ra C$. Since he proves
that for a Wirtinger cover the equality between Bs$(\bP\Gamma_{\tC })$
and $i(\tC)$ is scheme-theoretical outside the double points of
$i(\tC)$, it follows from his proof that, for a general double cover,
the base locus is the union of $i(\tC)$ and possibly a finite set of
points.  We prove (Sections 2 and 5)
\begin{thmi}\label{mainthm1}
If $g \geq 6$ and $C$ is non-hyperelliptic, the 
scheme-theoretical base locus in $P'$
of the linear system $\pp \Gamma_{\tC}$ is $i(\tC)$.
\end{thmi}

The proof of theorem 1 has two steps. First we show that
Bs$(\bP\Gamma_{\tC })$ equals $i(\tC)$ set-theoretically (Section 2).
In order to prove the scheme-theoretic equality, we introduce and
study divisors $D:= \Delta(E)$ in the linear systems $|2\Xi_0|$ and
$|2\Xi'_0|$ associated to certain semi-stable rank $2$ vector bundles
$E$ over the curve $C$ (Prop. 3.2). We calculate the tangent spaces to
the divisors $\Delta(E)$ along the curve $i(\tC)$ and show that at any
given point of $i(\tC)$ their intersection is equal to the tangent
space to $i(\tC)$.

\bigskip Let $\Theta_0$ be a symmetric theta divisor on the Jacobian
$JC$ and let $\alpha$ be the square-trivial invertible sheaf
associated to the double cover $\tC\ra C$. Translation by $\alpha$
induces an involution $T_\alpha$ on $JC$, which lifts canonically to a
linear involution acting on $H^0(JC,\T_0 + T^*_\alpha \T_0)$. Mumford
constructs in \cite{mum2} (see also \cite{vGP} Proposition 1)
canonical isomorphisms
\begin{equation} \label{mumiso}
\mu_+ : H^0(P,2\Xi_0) \map{\sim} H^0(JC, \T_0 + T^*_\alpha \T_0)_+ 
\qquad 
\mu_- : H^0(P,2\Xi_0) \map{\sim} H^0(JC, \T_0 + T^*_\alpha \T_0)_- 
\end{equation}
where the subscript $\pm$ denotes the $\pm$eigenspaces of the
involution. We are interested in some naturally defined subspaces of
these vector spaces. In connection with the Schottky problem, van
Geemen and van der Geer \cite{geem-geer} introduced the subspace
\[
\Gamma_{00} = \{ s\in H^0(A, 2\Theta) \: | \: \mult_0(s) 
\geq 4 \}
\]
for any abelian variety $A$ with symmetric principal polarization
$\Theta$.  It was conjectured by van Geemen, van der Geer and Donagi
(\cite{geem-geer} and \cite{donagi88} page 110) that if $(A,\Theta)$
is a Jacobian, then the base locus Bs$(\pp \Gamma_{00})$ of $\pp
\Gamma_{00}$ is the surface $C-C =\{\cO_C (p-q ) \ | \ p, q\in C\}
\subset JC$ as a set and, if $(A ,\Theta)$ is not in the closure of
the locus of Jacobians, then Bs$(\pp \Gamma_{00}) =\{ \cO \}$.
For Jacobians, the conjecture was proved by Welters \cite{welt1}.
For non-Jacobians, the conjecture was proved in dimension $4$ by the
first author \cite{izadi1}. Some evidence was also given for
non-Jacobian Pryms by the first author in \cite{izadi}. Consider
the subspaces $\Gamma_{\tilde{C}}^{(2)}$ of $H^0(P',2\Xi'_0)$ of
elements vanishing with multiplicity $\geq 2$ along $i(\tC)$ and the
subspace
\begin{equation*}
\Gamma^{\alpha}_{C-C} := \{ s\in  H^0 ( JC ,\T_0  + T^*_\alpha \T_0)
\: | \: C-C\subset Z(s)\} 
\end{equation*}
where $Z(s)$ denotes the zero divisor of the section $s$. This space
splits into $\pm$eigenspaces $\Gamma^{\alpha\pm}_{C-C}$ under the involution
induced by $T_{\alpha}$.

The infinitesimal study of the above mentioned divisors $\Delta(E)$
at the origin $\cO \in P$ and along the curve $i(\tC)$ allows us to prove
the following result (Section 4).
\begin{thmi}
  Assume $C$ non-hyperelliptic of genus $g\geq 5$. Via the canonical
  isomorphisms \eqref{mumiso}, we have equalities among the following
  subspaces
\begin{enumerate}
\item $\Gamma^{\alpha +}_{C-C} =\Gamma_{00}$, i.e., $\forall s \in
  H^0(P,2\Xi_0)$,
$$ \mult_0(s) \geq 4 \ \iff \ C-C \subset Z(\mu_+(s)) $$  
\item $\Gamma^{\alpha -}_{C-C} =\Gamma_{\tilde{C}}^{(2)}$, i.e., 
$\forall s \in H^0(P',2\Xi'_0)$,
$$ \Big( \forall \tp \in \tC \ \mult_{i(\tp)}(s) \geq 2 \Big)
 \ \iff \ C-C \subset Z(\mu_-(s)) $$  
\end{enumerate}
\end{thmi}
One can view statement $1$ as an analogue for Prym varieties of the
equivalence (see e.g. \cite{fay} page 489, \cite{welt1} prop 4.8 or
\cite{geem-geer})
$$ \forall s \in H^0(JC,2\T_0), \qquad \mult_0(s) \geq 4 \ \iff
\ C-C \subset Z(s). $$
Alternatively, one can derive equality $1$ from an analytic 
identity between Prym and Jacobian theta functions  (formula (41) \cite{fay}).
Equality $2$, however, seems to be new.
\bigskip

\noindent
{\em Acknowledgements:} We would like to thank R. Smith and R. Varley
for helpful discussions. This work was carried out during visits of
the first author of the University of Nice and of the second author of
the University of Georgia. We wish to thank these universities for
their hospitality. The first author was partially supported by a grant
from the National Security Agency.

\section{Preliminaries and notation}

In this section we introduce the notation and recall some well-known
facts on Prym varieties. Throughout the paper we will suppose the
genus of $C$ to be at least $5$. Let $\omega$ be the dualizing sheaf of $C$
and consider the two additional Prym varieties
$$ Nm^{-1}(\omega) = P_{even} \cup P_{odd} 
$$
which are characterized by the fact that 
$\dim H^0(\tilde{C},\lambda)$ is even (resp. odd) for $\lambda \in
P_{even}$ (resp. $P_{odd}$). The variety $P_{even}$ 
carries the naturally defined reduced
Riemann theta divisor
\[
\Xi :=\{\lambda\in P_{even}  \ | \  h^0 (\lambda ) > 0\}\; 
\]
a translate of which is $\Xi_0$.  Let $\modalpha$ and $\modkalpha$ be
the moduli spaces of semi-stable vector bundles of rank $2$ with
determinant $\alpha$ and $\omega\alpha$ respectively. Taking direct
image gives morphisms
$$ \varphi: P \cup P' \lra \modalpha \qquad \varphi: P_{even} \cup P_{odd} 
\lra \modkalpha. $$
Let
$\mathcal{L}_\alpha$ (resp. $\mathcal{L}_{\omega \alpha}$) be the generator
of the Picard group of $\modalpha$ (resp. $\modkalpha$). It is known that
$$
(\varphi_{|P})^* \mathcal{L}_\alpha = \cO(2\Xi_0) \qquad
(\varphi_{|P_{even}})^* \mathcal{L}_{\omega \alpha} = \cO(2\Xi).
$$
We denote by $\cO(2\Xi'_0)$ (resp. $\cO(2\Xi')$) the pull-back of
the line bundle $\mathcal{L}_\alpha$ (resp.  $\mathcal{L}_{\omega
  \alpha}$) to the Prym $P'$ (resp. $P_{odd}$), i.e.,
$(\varphi_{|P'})^* \mathcal{L}_\alpha = \cO(2\Xi'_0)$ and
$(\varphi_{|P_{odd}})^* \mathcal{L}_{\omega \alpha} = \cO(2\Xi')$.  We
consider the following morphisms
$$ \psi :  JC \lra \modalpha \qquad \xi \lms \xi \oplus \alpha \xi^{-1} $$
$$ \psi : \pic^{g-1}(C) \lra \modkalpha  \qquad \xi \lms \xi \oplus \omega 
\alpha \xi^{-1}. $$
One computes the pull-backs
$$ \psi^* \mathcal{L}_\alpha = \Theta_0 + T^*_\alpha \Theta_0 \qquad
\psi^* \mathcal{L}_{\omega \alpha} = \Theta + T^*_\alpha \Theta $$
where
\[
\Theta :=\{ L\in Pic^{ g-1} (C)  \ | \  h^0 (L )> 0\}\;
\]
and $\Theta_0$ is a symmetric theta divisor in the Jacobian $JC$,
i.e., a translate of $\Theta$ by a theta-characteristic. By abuse of
notation, we will also write $\mathcal{L}_{\alpha}$ and
$\mathcal{L}_{\omega \alpha}$ for $\psi^*\mathcal{L}_{\alpha}$ and
$\psi^*\mathcal{L}_{\omega\alpha}$ respectively. Note that $\psi$
induces linear isomorphisms at the level of global sections:
\begin{equation} \label{iso1}
\psi^* : H^0(\modalpha,\Lalpha) \cong H^0(JC,\Lalpha), \ \
\psi^* : H^0(\modkalpha,\Lomegaalpha) \cong H^0(\pic^{g-1}(C),\Lomegaalpha).
\end{equation} 
There is a well-defined morphism
$$ D : \modalpha \lra |\Lomegaalpha| \qquad (\text{resp.} \ 
 D : \modkalpha \lra |\Lalpha|) $$
where the support of $D(E)$ (reduced for $E$ general) is
$$ D(E) = \{ \xi \in JC \ (\text{resp.} \ \  \pic^{g-1}(C)) \: | \: 
h^0(C, E\otimes \xi) >0 \}. $$
\noindent
The two involutions of the Jacobian $JC$ given by 
$$ T_{\alpha}: \xi \lms \xi \otimes \alpha \qquad (-1): \xi \lms \xi^{-1} $$
induce (up to $\pm 1$) linear involutions $T_{\alpha}^*$ and $(-1)^*$ on the 
spaces of global sections
$H^0(JC,\Lalpha)$ and $H^0(\pic^{g-1}(C),\Lomegaalpha)$.
\begin{lem}
The projective linear involutions $T_{\alpha}^*$ and $(-1)^*$ acting on
$\pp H^0(JC,\Lalpha)$ are equal.
\end{lem}
\begin{proof}
  We observe that the composite map $T_{\alpha} \circ (-1) : \xi
  \mapsto \alpha \xi^{-1}$ verifies $\psi \circ (T_{\alpha} \circ
  (-1)) = \psi$. Since $\psi^*$ is a linear isomorphism \eqref{iso1},
  we have $(T_{\alpha} \circ (-1))^* = \pm id_{H^0}$.  Therefore
  $T_{\alpha}^* = \pm (-1)^*$.
\end{proof}
\noindent
Thus the two spaces decompose into $\pm$eigenspaces. Note that in
order to distinguish the two eigenspaces, we need a lift of the
$2$-torsion point $\alpha$ into the Mumford group. We will take the
following convention: the $+$eigenspace (resp. $-$eigenspace)
contains the Prym varieties $P$ and $P_{even}$ (resp. $P'$ and
$P_{odd}$), i.e., we have canonical (up to multiplication by a nonzero
scalar) isomorphisms.
\begin{equation} \label{iso2}
 H^0(JC,\Lalpha)_+ = H^0(P, 2\Xi_0), \qquad H^0(JC,\Lalpha)_- =
H^0(P',2\Xi'_0),  
\end{equation}
\begin{equation} \label{iso3}
 H^0(\pic^{g-1}(C),\Lomegaalpha)_+ = H^0(P_{even}, 2\Xi), \qquad
H^0(\pic^{g-1}(C), \Lomegaalpha)_- = H^0(P_{odd},2\Xi'). 
\end{equation}
\noindent
Since the surface $C-C$ is invariant under the involution $(-1) : \xi \mapsto
\xi^{-1}$, the subspace $\Gamma^{\alpha}_{C-C}$ is invariant under 
$(-1)^*$ and decomposes into a direct sum
of $\pm$eigenspaces for $(-1)^* = T_{\alpha}^*$:
$$
\Gamma^{\alpha}_{C-C} = \Gamma^{\alpha +}_{C-C} \oplus
\Gamma^{\alpha -}_{C-C}. $$

\bigskip
\noindent
{\em Prym-Wirtinger duality}
\bigskip

For the details see \cite{beau2} lemma 2.3. There exists an integral
Cartier divisor on the product $\modalpha \times \modkalpha$ whose
support is given by $$ \{ (E,F) \in \modalpha \times \modkalpha \: |
\: h^0(C, E \otimes F) >0
\}. $$
Its associated section can be viewed as an element of the tensor product
$$ H^0(\modalpha, \Lalpha) \otimes H^0(\modkalpha, \Lomegaalpha) $$
and it can be shown that the corresponding linear map 
\begin{equation} \label{wirtdual}
H^0(\modalpha, \Lalpha)^* \lra H^0(\modkalpha, \Lomegaalpha) 
\end{equation}
is an isomorphism and is equivariant for the linear involutions induced by the
map $E \mapsto E \otimes \alpha$. Hence 
using the identifications \eqref{iso2} and \eqref{iso3}
we obtain canonical isomorphisms,
\begin{equation} \label{wd1}
H^0(P, 2\Xi_0)^* \map{\sim}  H^0(P_{even}, 2\Xi) \qquad
H^0(P', 2\Xi'_0)^* \map{\sim}  H^0(P_{odd}, 2\Xi'). 
\end{equation}

\section{The base locus of $\pp \Gamma_{\tC}$}

In this section we compute the set-theoretical base locus of the
subseries $\pp \Gamma_{\tC}$ on the Prym variety $P'$. Suppose $C$
non-hyperelliptic. We need some additional notation:
denote by $\tC_m$ the $m$-th symmetric power of $\tC$ and let $S$ be
the subvariety of $\tC_{2g-2}$ defined as
$$ S = \{ D \in \tC_{2g-2} \ | \ Nm(D) \in |\omega| \ \text{and} \
h^0(D) \equiv 1 \ \text{mod} \ 2 \}. $$
Then, by \cite{beau} Corollaire page 365, the variety $S$ is normal and
irreducible of dimension $g-1$. The variety $S$ comes equipped with
two natural surjective morphisms
$$ Nm : S \lra |\omega| \qquad u : S \lra P_{odd} $$
where $u$ associates to an effective divisor $D$ its line
bundle $\cO_{\tC}(D)$. Note that $u$ is
birational and $Nm$ is finite of degree $2^{2g-3}$. 
Also denote by $u$ the extended morphism $u : \tC_{2g-2} \ra
\pic^{2g-2}(\tC)$ and consider the commutative diagram
\begin{equation} \label{comdiags}
\begin{array}{ccc}
S & \inj & \tC_{2g-2} \\
 & & \\
\downarrow^u & & \downarrow^u \\
 & & \\
P_{odd} & \inj & \pic^{2g-2}(\tC).
\end{array}
\end{equation}
Consider the Brill-Noether locus in $P_{odd}$ which is
defined set-theoretically by
$$ \Xi_3 := \{ \lambda \in P_{odd} \ | \ h^0(\lambda) \geq 3 \}. $$
The scheme structure on $\Xi_3$ is defined by taking the
scheme-theoretical intersection \cite{welt2}
$$
\Xi_3 := W^2_{2g-2}(\tC) \cap P_{odd} $$
where $W^2_{2g-2}(\tC)
\subset Pic^{2g-2}\tC$ is the Brill-Noether locus of line bundles
having at least $3$ sections (see \cite{acgh}).
\begin{lem}
The subscheme $\Xi_3 \subset P_{odd}$ is not empty and is
of pure codimension $3$.
\end{lem}

\begin{proof}
Theorem 9 \cite{dcp} asserts that $\Xi_3$ is not empty and every
irreducible component has dimension at least $g-4$. Suppose that
there is an irreducible component $I$ of dimension $\geq g-3$. Then
its inverse image $u^{-1}(I)$ has dimension $\geq g-1$, hence, since
$S$ is irreducible, $u^{-1}(I) = S$ and $\Xi_3 = P_{odd}$. The
last equality can not happen, since otherwise, using translation by an
element of the form $\cO_{\tC}(\tp - \sigma \tp)$, we would have
$\Xi = P_{even}$.
\end{proof}

Observe that $u$ is equivariant for the action of $\sigma$ on $S$ and
$P_{odd}$. Denote by $Z = u^{-1}(\Xi_3)$ the inverse image of the
subscheme $\Xi_3$. By the previous lemma $Z$ is of pure codimension
$1$ in $S$. We will see in a moment that there is a Cartier divisor
$\cD$ on $S$ whose support is the support of $Z$.  Let $\omega_{\tC}$
be the dualizing sheaf of $\tC$. Consider the following divisors in
$\tC_{2g-2}$
\[
U_{\tp} := \{ D \in \tC_{2g-2} \ | \ \exists D'\in \tC_{ 2g-3 } \text{
  with } D = D'+\tp\}
= \tp + \tC_{2g-3}
\]
\[
V_{\tp} := \{ D \in \tC_{2g-2} \ | \ h^0(\omega_{\tC}(-D-\tp)) \geq 1
\}
\]  
and let $\bar{U}_{\tp}$ and $\bar{V}_{\tp}$ be their intersections 
with $S$. A straightforward calculation involving Zariski
tangent spaces then shows that $\bar{U}_{\tp}$ is a reduced 
divisor.
We denote by $\cO_S(1)$ the pull-back by the norm map of the hyperplane
line bundle on $|\omega|$. Then it is easily seen that, for any 
$\tp \in \tC$,
\begin{equation} \label{pbnm}
Nm^*(|\omega(-p)|) = \bar{U}_{\tp} + \bar{U}_{\sigma \tp} \in
|\cO_S(1)|.
\end{equation}
Let $\tTheta_{\lambda}$ denote the translate of
$\tTheta$ by $\lambda$. Then, for any points $\tp, \tq \in \tC$, we
have an equality among divisors on $\tC_{2g-2}$ (see \cite{welt1} page 6)
\begin{equation} \label{pbthetapq}
u^*(\tTheta_{\tp - \tq}) = U_{\tp} + V_{\tq}.
\end{equation}

The analogue on the even Prym variety of the following lemma was
previously proved by R. Smith and R. Varley. In the case of genus $3$
it is in their paper \cite{smithvarley99} (Prop. 1 page 358) and for
higher genus it will be published in their upcoming paper
\cite{smithvarley00}.
\begin{lem}
There exists an effective Cartier divisor $\cD$ on $S$ whose 
support is equal to
$$ \text{supp} \ Z = \{ D \in S \ | \ h^0(\cO_{\tC}(D)) \geq 3 \}. $$
Moreover, we have the following equality among effective
Cartier divisors 
\begin{equation} \label{eqxiD}
u^* (\Xi_{\tp -\sigma\tp} + \Xi_{\sigma\tp -\tp} ) =
\bar{U}_{\tp} + \bar{U}_{\sigma \tp} + \cD \ \ \ \forall \tp \in \tC.
\end{equation}
In particular, $u^* \cO_{P_{odd}}(2\Xi') = \cO_S(1) \otimes
\cO_S(\cD)$.
\end{lem}

\begin{proof}
  We are going to define $\cD$ as the residual divisor of the
  restricted divisor $\bar{V}_{\tq}$, for a given point $\tq \in \tC$
  and then show that it does not depend on the choice of $\tq$. We
  first observe that we have an equality of sets
$$\bar{V}_{\tq} = \bar{U}_{\sigma \tq} \cup Z $$
which can be seen as follows: for $D\in\tC_{ 2g-2 }$
such that $h^0 (D) = h^0(\omega_{\tC}(-D)) = 1$ 
the
assumption $D \in \bar{V}_{\tq}$ and the formula $D + \sigma D
= \pi^* (Nm (D))$ imply that
$ \tq \in \text{supp} \ \sigma D \ \iff \
D \in \bar{U}_{\sigma\tq}$ .
If $h^0(D) = h^0(\omega_{\tC}(-D)) \geq 2$, then $D \in \text{supp}\ Z$.
Again a calculation involving Zariski tangent spaces shows that
$\bar{V}_{\tq}$ is reduced generically on $\bar{U}_{\sigma \tp}$.
Hence we can define
$\cD$ by $\bar{V}_{\tq} = \bar{U}_{\sigma \tq} + \cD$. Now
we substitute this expression into \eqref{pbthetapq}, which we
restrict to $S$
$$
u^*(\tTheta_{\tp - \tq})_{|S} = \bar{U}_{\tp} + \bar{U}_{\sigma
  \tq} + \cD.$$
Now we fix $\tq$ and we take the limit when $\tp \ra
\tq$. Since $\cO_{P_{odd}}(\tTheta) = \cO_{P_{odd}}(2\Xi')$, we see
that $\bar{U}_{\tp} + \bar{U}_{\sigma \tp} + \cD \in |u^*\cO_{P_{odd}}
(2\Xi')|$. So by \eqref{pbnm} we get the line bundle equality claimed
in the lemma and we see that the scheme-structure on $\cD$ does not
depend on the point $\tq$. To prove \eqref{eqxiD}, we compute using
\eqref{pbthetapq}
$$u^* (\tTheta_{\tp -\sigma\tp} + \tTheta_{\sigma\tp -\tp}) =
\bar{U}_{\tp} + \bar{V}_{\sigma \tp} + \bar{U}_{\sigma \tp} +
\bar{V}_{\tp}.$$
Now we restrict to $S$ and use the commutativity of
diagram \eqref{comdiags} and the divisorial equality $\tTheta_{\tp -
  \sigma \tp} \cap P_{odd} = 2 \Xi_{\tp - \sigma \tp}$ to obtain
$$u^* (2\Xi_{\tp -\sigma\tp} +  2\Xi_{\sigma\tp -\tp}) =
2\bar{U}_{\tp} + 2\bar{U}_{\sigma \tp} + 2\cD.$$
Since $\bar{U}_{\tp} + \bar{U}_{\sigma \tp} + \cD \in |u^*\cO_{P_{odd}}
(2\Xi')|$ we can divide this equality by $2$ and we are done.  
\end{proof}

Let $\mu$ be a point of $Bs(\pp \Gamma_{\tC})$. By lemma \ref{lemtaui}
the linear map $i^* : |\omega|^* \lra |2 \Xi'_0|^*$ is injective and,
since $|\omega|^*$ is the span of the image of $\tC$ in $|2\Xi'_0|^*$,
the space $\pp \Gamma_{\tC}$ is the annihilator of $|\omega|^* \subset
|2\Xi'_0|^*$. So $Bs(\pp \Gamma_{\tC}) = |\omega |^*\cap Kum(P')$ and
$\mu$ corresponds to a hyperplane $H_\mu \in |\omega|^*$. Since
$\mu\in Kum(P')$, the image of $\mu$ by Wirtinger duality is the
divisor $\Xi_{\mu} +\Xi_{\mu^{ -1}}\in |2\Xi' |$.

\begin{lem}
With the previous notation, we have an equality
\begin{equation} \label{pbhmu}
\forall \mu \in Bs(\pp \Gamma_{\tC}) \qquad 
Nm^*(H_\mu) + \cD = u^*(\Xi_\mu + \Xi_{\mu^{-1}}).
\end{equation}
\end{lem}

\begin{proof} 
  The equality follows from the commutativity of the right-hand square
  of the diagram
$$
\begin{array}{ccccccccc}
\tC & \map{\pi} & C & \map{\varphi_{can}} & |\omega|^* &
\map{Nm^*} & |\cO_S(1)| & \map{+\cD} & |\cO_S(1) \otimes
\cO_S(\cD)| \\
 & & & & & & & & \\
\downarrow^i & & & & \downarrow^{i^*} & & & \nearrow^{u^*} &  \\
 & & & & & & & &  \\
P' & & \lra  &   & |2\Xi'_0|^* & \cong & |2\Xi'|. & & 
\end{array}
$$
The commutativity of the right-hand square follows from that of the
outside square because $\varphi_{can}(C)$ generates $|\omega|^*$. In
other words we need to check the assertion of the lemma only for
hyperplanes of the form $|\omega(-p)|$ for $p \in C$. This follows
immediately from \eqref{pbnm} and \eqref{eqxiD}.
\end{proof}

\begin{cor}
For every $\mu\in Bs(\pp \Gamma_{\tC})$, the hyperplane $Nm^* (H_{\mu
  })$ is reducible.
\end{cor}
\begin{proof} By the above Lemma we have
\[
u^*(\Xi_\mu + \Xi_{\mu^{-1}}) -\cD = Nm^*(H_\mu).
\]
If $Nm^*(H_\mu)$ is irreducible, then the support of one of the
divisors $u^*(\Xi_\mu )$ or $u^*(\Xi_{\mu^{-1}})$, say $u^*(\Xi_\mu
)$, is contained in the support of $\cD$. This is impossible because
$u^*(\Xi_\mu )$ is the inverse image of a divisor in $P_{ odd}$ and
$supp \cD$ is the inverse image of the codimension $3$ support of
$\Xi_3$.
\end{proof}

The set-theoretical assertion of
theorem \ref{mainthm1} now follows from the following lemma.

\begin{lem}\label{lemHred}
If $C$ is not bi-elliptic, we have a set-theoretical equality
$$
\{ H \in |\omega|^* \ : \ Nm^*(H) \ \text{is reducible} \} =
\varphi_{can}(C). $$
If $C$ is bi-elliptic, the LHS is contained in the
union of $\varphi_{can}(C)$ and the finite set of points $t\in |\omega|^*$
such that the projection from $t$ induces a morphism of degree $2$
from $C$ onto an elliptic curve.
\end{lem}
\begin{proof}
  Suppose that $Nm^*(H)$ is reducible. Then a local computation shows
  that the hyperplane $H$ is everywhere tangent to the branch locus of
  $Nm$. It is immediately seen that the branch locus $B$ of $Nm$ is
  the dual hypersurface of the canonical curve. The components of the
  singular locus $Sing(B)$ of $B$ are of two different types which can be
  described as follows
\begin{description}
\item[type 1] whose points are hyperplanes tangent to 
  $\varphi_{can}(C)$ in more than one point.
\item[type 2] whose points are hyperplanes osculating to
  $\varphi_{can}(C)$.
\end{description}
To prove that $\mu\in\varphi_{can}(C)$, we need to prove that there is
a point on $H \cap B$ which is smooth on $B$ because the dual
variety of $B$ is the closure of the set of hyperplanes tangent to $B$
at a smooth point and this is equal to $\varphi_{can}(C)$. In other
words we need to show that $H \cap B$ is not contained in
$Sing (B)$. Since $H\cap B$ has pure codimension $2$, it suffices to show
that no codimension $2$ component of $Sing (B)$ is contained in a hyperplane.

Suppose a codimension $2$ component $B_i$ of type $i$ ($i=1$ or $2$)
is contained in a hyperplane $H$ in $|\omega |$ and let $t\in |\omega
|^*$ be the corresponding point. Then the set of hyperplanes in
$|\omega |^*$ through $t$ and doubly tangent (resp. osculating) to
$\varphi_{ can }(C)$ has dimension $g-3$. We have
\begin{lem}
  For any $t\in\varphi_{ can }(C)$ the restriction $\rho$ of the
  projection from $t$ to $\varphi_{ can }(C)$ is birational onto its
  image. If $t\in |\omega |^*\setminus\varphi_{ can }(C)$, then $\rho$
  is either birational onto its image or of degree two onto an
  elliptic curve.
\end{lem}
\begin{proof} 
  First note that the degree of the image $C_t$ of $C$ by the
  projection is at least $g-2$ because $C_t$ is a non-degenerate curve
  in a projective space of dimension $g-2$. If $t\in\varphi_{ can
    }(C)$, then the degree of $\rho$ is equal to $2g-3$. The degree
  $r$ of the restriction of $\rho$ to $C_t$ verifies
  $r\cdot\hbox{deg}(C_t) = 2g-3$. Therefore $\frac{2g-3 }{r}\geq g-2$.
  Or $r\leq 2 +\frac{1 }{g-2 }$ which implies $r\leq 2$. However, $r$
  cannot be equal to $2$ because $2g-3$ is odd. If
  $t\not\in\varphi_{ can }(C)$, then the same argument gives again
  $r\leq 2$ because $g\geq 5$. Hence, if $\rho$ is not generically injective,
  then $r =2$ and $\hbox{deg} (C_t ) = g-1$. Therefore $C_t$ is either 
  smooth rational or an elliptic curve. Since $C$ is not hyperelliptic, 
  we have that $C_t$ is an elliptic curve.
\end{proof}

First suppose that $C\ra C_t$ is birational. If $i=1$, projecting from
$t$, we see that the set of hyperplanes in $|\omega |^*/t$ doubly
tangent to $C_t$ has dimension $(g-3)$=dimension of the dual variety of
$C_t$ which is impossible. If $i=2$, then the set of hyperplanes in
$|\omega |^*/t$ osculating to $C_t$ has dimension $g-3$ which is also
impossible.

If $C\ra C_t$ is of degree $2$, then indeed every hyperplane tangent
to $C_t$ pulls back to a hyperplane twice tangent (or osculating if
the point of tangency is a branch point of $C\ra C_t$) to $\varphi_{
  can } (C)$ and we have a codimension $2$ family of type $B_1$
contained in the hyperplane corresponding $H$ to $t$. Then $Nm^* (H)$
could be reducible.
\end{proof}

The previous lemma proves theorem 1 set-theoretically
 for a non bi-elliptic curve.  In
the bi-elliptic case, we have to work a little more. By Lemma
\ref{lemHred} a hyperplane $H \not\in \varphi_{can}(C)$, such that
$Nm^*(H)$ might be reducible, corresponds to a point $e \in
|\omega|^*$ such that the projection from $e$ induces a morphism
$\gamma$ of degree $2$ from $C$ to an elliptic curve $E$. In other
words, $e$ is the common point of all chords $\langle \gamma^* q
\rangle ; \ (q \in
E)$. In that case there exists a 
1-dimensional family (parametrized by $E$) of trisecants,
namely the chords $\langle \gamma^*q \rangle$, to the Kummer variety
$Kum(P')$. By \cite{debarre} the Prym variety is a Jacobian and by
\cite{sho} (see also \cite{beau3} page 610) the double cover 
$\pi : \tC \ra C$ is of the following two types
\begin{enumerate}
\item $C$ is trigonal
\item $C$ is a smooth plane quintic and $h^0(\cO_C(1)\otimes \alpha) = 0$
\end{enumerate}

\begin{lem}
No double cover of a  bi-elliptic curve $C$ of genus $g \geq 6$ 
is of the above two types. 
\end{lem}

\begin{proof}
  For a bi-elliptic curve $C$, the Brill-Noether locus $W^1_{g-1}(C)$
  has two irreducible components, which are fixed by the reflection in
  $\omega$ (\cite{welt1} Corollary 3.10). For a smooth plane quintic
  this Brill-Noether locus is irreducible, ruling out 1.  For a
  trigonal curve this Brill-Noether locus has two irreducible
  components, which are interchanged by reflection in $\omega$, ruling
  out 2.
\end{proof}

\begin{rem}
  If $g = 5$ and $C$ is bi-elliptic, we do not know whether the common
  point of all the chords for a given bi-elliptic structure lies on
  $Kum(P')$ (see also \cite{beau3} Remark (1) page 611). We expect it
  not to be on $Kum(P')$.
\end{rem}

\section{Rank $2$ bundles and $2\Xi$-divisors}

Consider the induced action of the involution $\sigma$ on the 
moduli space
$\mathcal{SU}_{\tC}(2,\cO)$ given by $\tE\mapsto \sigma^* \tE$.
Since the covering $\pi$ is unramified, the fixed point set for
the $\sigma$-action
$$Fix_{\sigma} \mathcal{SU}_{\tC}(2,\cO) = \{ [\tE] \in
\mathcal{SU}_{\tC}(2,\cO) \ | \ \exists \theta : \ \sigma^* \tE \map{\sim}
\tE \} $$  
has two connected components which are the isomorphic images of
$\mod0$ and $\modalpha$ by $\pi^*$. Similarly, since 
$\sigma^* \omega_{\tC} \map{\sim} \omega_{\tC}$, the
involution $\sigma$ acts on
$\mathcal{SU}_{\tC}(2,\omega_{\tC})$ and
$$Fix_{\sigma} \mathcal{SU}_{\tC}(2,\omega_{\tC}) = \pi^*\modk \cup
\pi^*\modkalpha $$

\begin{prop}\label{propDelta}
Consider a bundle $E \in \modkalpha$ such that $E \not\in \varphi(P_{odd})$
and put $\tE = \pi^* E$. Then there is a divisor $\Delta(E) \in
|2\Xi_0|$ with the following properties.
\begin{enumerate}
\item If $D(\tE)$ does not contain $P$, then 
\[
D(\tE) = 2\Delta(E).
\]
For $E$ general, $P$ is not
contained in $D(\tE)$ and $\Delta(E)$ is reduced.

\item Let $pr_+$ be the projection $|\Lalpha|\ra |2\Xi_0|$ with center 
  $| 2\Xi_0 '|$ (see (\ref{iso2})). Then we
  have a commutative diagram
$$
\begin{array}{ccc}
\modkalpha & \map{D} & |\Lalpha| \\
 & & \\
 & \searrow^{\Delta} & \downarrow^{pr_+} \\
 & & \\
 & & |2\Xi_0| = |\Lalpha|_+
\end{array}
$$ 
\end{enumerate}
\end{prop} 

\begin{rem}
Similarly, when  $E\in\modkalpha$ such that $E\not\in\varphi(P_{even})$,
we get divisors $\Delta(E) \in |2\Xi'_0|$ as described in
prop. 3.2 by
projecting on the $-$eigenspace $pr_- : |\Lalpha| \lra
|\Lalpha|_- = |2\Xi'_0|$.
\end{rem}

\begin{proof}
  1. Given a bundle $F \in Fix_\sigma
  \mathcal{SU}_{\tilde{C}}(2,\omega_{\tilde{C}})$ and a line bundle
  $\xi \in J\tilde{C}$ which is anti-invariant under $\sigma$, i.e.,
  $\sigma^* \xi \map{\sim} \xi^{-1}$, we have a natural non-degenerate
  quadratic form with values in the canonical bundle
  $\omega_{\tilde{C}}$
\begin{eqnarray*}
q: F \otimes \xi & \lra & \omega_{\tilde{C}} \\
 s & \lms & s \wedge \sigma^* s
\end{eqnarray*}
where $s$ is a local section of $F \otimes \xi$. Note that we have
canonical isomorphisms
$$ \sigma^* (F\otimes \xi) = F \otimes \xi^{-1} = 
\mathrm{Hom}(F\otimes \xi,\omega_{\tilde{C}}) $$
Therefore we are in a position to apply the Atiyah-Mumford lemma 
\cite{mum} to
the family of bundles (here $F$ is fixed, with $\sigma^* F \map{\sim} F$)
$$  \{F \otimes \xi \}_{\xi \in P }$$
which states that the parity of $h^0(\tilde{C}, F\otimes \xi)$ is
constant when $\xi$ varies in $P$.

\bigskip
\noindent
>From now on, we suppose $F=\tE = \pi^*E$, with $E \in \modkalpha$, then 
$$
h^0(\tilde{C}, \tE) = 2h^0(C,E) \equiv 0 \ \ \text{mod} \ \ 2. $$
For the first equality we use the fact that $H^0(\tilde{C},\tE) =
H^0(C,E) \oplus H^0(C, E\alpha)$ and, by Riemann-Roch and Serre duality,
$h^0(C,E) = h^1 (C, E) = h^0 (C ,\omega\otimes E^* ) = h^0 (C, E\alpha
)$.

First suppose that $E\in\modkalpha$ is general. Then the divisor
$D(\tE)$ does not contain the Prym variety $P$ (e.g. because, for
general $E$, $h^0(E)=0 \ \iff \ h^0 (\tE )=0 \ \iff \ \cO \not\in
D(\tE)$), so the restriction of the divisor $D(\tE) \in |2
\Theta_{\tC}|$ to $P$ is a divisor in the linear system $|4\Xi_0|$.
Moreover, for $\xi \in D(\tE)\cap P$
$$\mult_\xi D(\tE) \geq h^0(\tilde{C}, \tE \otimes \xi) \geq 2 $$
because $h^0(\tilde{C}, \tE \otimes \xi) \equiv h^0(\tilde{C}, \tE )
\equiv 0 \ \ \text{mod} \ \ 2$.  Hence any point $\xi \in D(\tE)\cap
P$ is a singular point of $D(\tE)$, which implies that $D(\tE) \cap P$
is an everywhere non-reduced divisor. We have
\begin{lem}
  Suppose that $D(\tE) \cap P$ is a divisor in $P$. Then there is a
  divisor $\Delta(E) \in |2\Xi_0|$ such that $D(\tE) \cap P = 2
  \Delta(E)$.
\end{lem}
\begin{proof}
  A local equation of $\Delta(E)$ is given by the pfaffian of a
  skew-symmetric perfect complex of length one $L \lra L^*$
  representing the perfect complex $Rpr_{1*}(\mathcal{P} \otimes
  pr^*_2 \tE)$ where $\mathcal{P}$ is the Poincar\'e line bundle over
  the product $P \times \tC$ and $pr_1, pr_2$ are the projections on
  the two factors. The construction of the complex $L \lra L^*$ is
  given in the proof of Proposition 7.9 \cite{ls}.
\end{proof}
If $E$ is of the form $E = \pi_* L$ for some $L \in P_{even}$,
we have $\Delta(E) = T^*_L \Xi + T^*_{\omega L^{-1}} \Xi$. It 
follows from this equality that $\Delta(E)$ is reduced for
general $E$.

So far we have defined a rational map $\Delta: \modkalpha \lra
|2\Xi_0|$. It will follow from part 2 of the proposition that $\Delta$
can be defined away form $\varphi(P_{odd})$.

\bigskip
\noindent
2. First we consider the composite (rational) map
$$
\pic^{g-1}(C) \map{\psi} \modkalpha \map{\Delta} |2\Xi_0|. $$
A
straight-forward computation shows that for all $\xi \in
\pic^{g-1}(C)$ such that $\pi^*\xi\not\in P_{odd}$ the divisor
$\Delta(\psi(\xi)) = \Delta(\xi \oplus \omega \alpha \xi^{-1})$ equals
the translated divisor $T_{\pi^* \xi} \tilde{\Theta}$ restricted to
$P$. Hence, by \cite{mum2}, the map $\Delta \circ
\psi$ is given by the full linear system $|\Lomegaalpha|_+$ of
invariant elements of $|\Lomegaalpha|$. By Prym-Wirtinger duality
\eqref{wirtdual} and \eqref{wd1} $|\Lomegaalpha|^*_+ \cong |\Lalpha|_+
\cong |2\Xi_0|$ and we obtain the commutative diagram in the proposition.
Geometrically, $\Delta$ is obtained by restricting the projection with
center the $-$eigenspace $|\Lalpha|_-$ to the embedded moduli space
$\modkalpha \subset |\Lalpha|$.  Since by \cite{nr} $|\Lalpha|_- \cap
\modkalpha = \varphi(P_{odd})$ we see that $\Delta$ is well-defined
for $E\not\in\varphi(P_{odd})$ even if $D(\tE)\supset P$.
\end{proof}

\begin{rem}
We observe that we obtain by the same construction a rational map
$$ \Delta : \modk \lra |2 \Xi_0|\; . $$
The images under $\Delta$ of the two moduli spaces $\modk$ and $\modkalpha$
coincide, which is easily deduced from the following formula. Let 
$\beta$ be a $4$-torsion point such that $\beta^{\otimes 2} = \alpha$ 
and $\pi^* \beta
\in P[2]$. Then, for any $E \in \modk$, we have $E \otimes \beta \in
\modkalpha$ and
$$ T_{\pi^* \beta}^*  \Delta(E) = \Delta(E\otimes \beta)\; . $$
\end{rem}

\noindent
Similar statements hold for $\modalpha$.

\section{Proof of theorem 2}

\subsection{Proof of $\Gamma^{\alpha +}_{C-C} = \Gamma_{00}$}
The strategy is to show that the two linear maps 
$$\phi_1 :  H^0(P,2\Xi_0)_0 \lra \sym^2 T_0^*P = \sym^2 H^0(\omega \alpha)$$
and
$$\phi_2 : H^0(JC,\Lalpha)_{+0} \lra H^0(C\times C, \delta^* \Lalpha -
2\Delta)_+ = \sym^2 H^0(\omega \alpha)$$
differ by multiplication by a scalar under the
isomorphism \eqref{iso2} $H^0(JC,\Lalpha)_{+0 }\cong H^0(P,2\Xi_0)_0$.
Here the subscript $0$ denotes the subspace (on $P$ or $JC$)
consisting of global sections vanishing at the origin. The map
$\phi_1$ sends $s\in H^0(P,2\Xi_0)_0$ to the quadratic term of its
Taylor expansion at the origin $\cO \in P$ and $\phi_2$ is the
pull-back of invariant sections of $\Lalpha$ under the difference map
\begin{eqnarray*}
\delta : C \times C & \lra & JC \\
(p,q) & \lms & \cO_{C}(p-q).
\end{eqnarray*}
By restricting to the fibers of the two projections $p_i :C\times C\ra
C$ and using the See-saw Theorem, we compute that $\delta^* \Lalpha =
p_1^* (\omega \alpha) \otimes p_2^* (\omega \alpha )(2 \Delta_C)$
where $\D_C\subset C\times C$ is the diagonal. Since $\phi_2^{ -1} (0)
=\D_C$ and the sections of $\cL_{\alpha }$ are symmetric, we see that
$\im \phi_2 \subset \sym^2 H^0(\omega \alpha )\subset H^0
(\omega\alpha )^{\otimes 2} = H^0 (p_1^* (\omega \alpha) \otimes p_2^*
(\omega \alpha ))\subset H^0 (p_1^* (\omega \alpha) \otimes p_2^*
(\omega \alpha ) (2\D_C ))$. So if $\phi_1$ and $\phi_2$ are
proportional, we will have
$$ \Gamma_{00} = \ker \phi_1 = \ker \phi_2 = \Gamma^{\alpha +}_{C-C}. $$
\noindent
To show that $\phi_1 = \lambda\phi_2$ for some $\lambda\in\bC^*$, we
compute $\phi_1(s_E)$ and $\phi_2(s_E)$ for special sections, namely those
with divisor of zeros $Z(s_E) = \Delta(E)$ for some
vector bundle $E \in \modkalpha$ with $h^0(E) = h^0(E
\otimes \alpha) = 2$. Recall that by Riemann-Roch and Serre
duality we have for $h^0(E) = h^0(E \otimes
\alpha)$ for $E \in \modkalpha$. Now to compute $\phi_1(s_E)$, 
we need to determine the
tangent cone to $\Delta(E)$ at $\cO \in P$.  As before we put
$\tE = \pi^*E$. By \cite{laszlo} prop. V.2, this tangent cone is the
intersection of the anti-invariant part $H^0(\omega \alpha) =
H^0(\omega_{\tilde{C}})_-$ of $H^0(\omega_{\tC}) = T^*_0 J\tilde{C}$
with the affine cone over the projective cone over 
the Grassmannian $Gr(2, H^0(\tE)^*) \subset \pp \Lambda^2
H^0(\tE)^*$ under the linear map
\begin{equation} \label{mapmu}
\mu^* : H^0(\omega_{\tilde{C}})^* \lra \Lambda^2 H^0(\tE)^*
\end{equation}
which is the dual of the map $\mu :\Lambda^2 H^0(\tE)\ra
H^0(\omega_{\tilde{C}})$ obtained from exterior product by the
isomorphism $\wedge^2\tE\cong\omega_{\tC }$. Note that the
$\sigma$-invariant part $[\Lambda^2 H^0(\tE)^*]_+$ is 
canonically isomorphic to the $2$-dimensional subspace
$\Lambda^2 H^0(E)^* \oplus \Lambda^2 H^0(E\alpha)^* \subset
\Lambda^2 H^0(\tE)^*$ because $H^0(\tE)_+ = H^0(E)$ and $H^0(\tE)_- =
H^0(E \alpha)$. Since $\wedge^2 E\cong\wedge^2 (E\otimes\alpha
)\cong\omega\alpha$, the restriction of $\mu$ to $\wedge^2 H^0 (E)$
(resp. $\wedge^2 H^0 (E\otimes\alpha )$) which is obtained from
exterior product by the isomorphism $\wedge^2 E\cong\omega\alpha$
(resp. $\wedge^2 (E\otimes\alpha )\cong\omega\alpha$) maps into $H^0
(\omega\alpha )$. Therefore
the linear map $\mu^*$ \eqref{mapmu} maps $\sigma$-anti-invariant
sections into $\sigma$-invariant sections, i.e.,
\begin{equation} \label{mapmuplus}
\mu^*_+ : H^0(\omega \alpha)^* \lra \Lambda^2 H^0(E)^* \oplus \Lambda^2 
H^0(E\alpha)^*.
\end{equation}
Since the intersection $\pp (\Lambda^2 H^0(E)^* \oplus \Lambda^2
H^0(E\alpha)^*) \cap Gr(2,H^0(\tE)^*)$ consists of the two points $\pp
\Lambda^2 H^0(E)^*$ and $\pp \Lambda^2 H^0(E\alpha)^*$, it follows
that the intersection of $H^0 (\omega\alpha )\subset H^0 (\omega_{\tC
  })$ with the cone over $Gr(2, H^0(\tE)^*)$ is the union of the two
lines $\wedge^2 H^0 (E)$ and $\wedge^2 H^0 (E\otimes\alpha )$.
Therefore the tangent cone of $\Delta(E)$ at the origin is the union
of the two hyperplanes in $|\omega \alpha|^*$ which are the zeros of
$a,b \in H^0(\omega \alpha)$ such that
\begin{equation} \label{defab}
a \cc = \im (\Lambda^2 H^0(E) \lra H^0(\omega \alpha)) \qquad
b \cc = \im (\Lambda^2 H^0(E\alpha) \lra H^0(\omega \alpha)).
\end{equation}
In other words, up to multiplication by a nonzero scalar,
$$ \phi_1(s_E) = a \otimes b + b \otimes a \in \sym^2 H^0(\omega
\alpha). $$
\noindent
We now compute $\phi_2(s_E)$.  First we note that the pull-back map
induced by $\delta$ is equivariant for the involution $(-1) : \xi
\mapsto \xi^{-1}$ acting on $JC$ and the involution $(p,q) \mapsto
(q,p)$ acting on $C \times C$. Since $\D(E) = pr_+ (D(E))$ by
Proposition \ref{propDelta}, this implies that
\begin{equation} \label{comprdelta}
\phi_2 (s_E) =\phi_2(pr_+(s_E)) = pr_+ (\delta^* (s_E)) 
\end{equation}
On the RHS $pr_+$ denotes the projection $H^0(\omega \alpha) \otimes
H^0(\omega \alpha) \lra \sym^2 H^0(\omega \alpha)$. Therefore we
compute
$$\delta^* (D(E)) = \{ (p,q) \in C \times C \: | \:  h^0(E(p-q))>0 \} $$
and take its symmetric part. It follows from \cite{geem-iz} lemma 3.2
that
\begin{equation} \label{pbdelta} 
\delta^* (D(E)) =  C \times Z_a + Z_b \times C + 2 \Delta_C
\end{equation}
where $Z_a$ (resp. $Z_b$) is the divisor of zeros of $a$ (resp. $b$).
Hence it follows from \eqref{comprdelta} and \eqref{pbdelta} that
$\phi_2(pr_+(s_E)) = a \otimes b + b \otimes a$ up to multiplication
by a nonzero scalar. We can now conclude that $\phi_1 = \lambda\phi_2$
for some
$\lambda\in\bC^*$ because, by the following lemma (prop. 3.7
\cite{geem-iz}), we have enough bundles $E \in \modkalpha$ with
$h^0(E) = 2$ to generate linearly the image $\sym^2 H^0(\omega
\alpha)$ of $\phi_1$ and $\phi_2$.

\begin{lem}(prop. 3.7 \cite{geem-iz})
  For general sections $a,b \in H^0(\omega \alpha)$, we can find a
  semi-stable bundle $E \in \modkalpha$ with $h^0(E)=2$ such that
  \eqref{pbdelta} holds.
\end{lem}

\subsection{Proof of $\Gamma^{\alpha -}_{C-C} = \Gamma^{(2)}_{\tC}$}
First note that any anti-invariant section of
$\Lalpha$ vanishes at $\cO \in JC$. Denote by 
$$
\tau : H^0(JC,\Lalpha)_- \lra T_0^*JC = H^0(\omega) $$
the map
which sends an element $s$ of $H^0(JC,\Lalpha)_-$ to the linear term
of its Taylor expansion at the origin (Gauss map). Recall the
natural embedding of the curve $\tilde{C}$ into the Prym variety $P'$
\begin{equation} \label{embtc}
i: \tilde{C}  \lra  P' \qquad
 \tilde{p}  \lms  \cO_{\tC}(\tilde{p} - \sigma \tilde{p}).
\end{equation}  
Then $i^* \cO(2\Xi'_0)\cong \omega_{\tilde{C}}$ and since all
$2\Xi'_0$-divisors are symmetric and $i$ is equivariant for
the involution, $i$ induces a
linear map
\begin{equation} \label{pbi}
i^* : H^0(P',2 \Xi'_0) \lra H^0(C, \omega) = H^0(\tC, \omega_{\tC})_+
\end{equation}

\begin{lem}\label{lemtaui}
The linear maps $\tau$ and $i^*$ are proportional via the isomorphism
\eqref{iso2} and are surjective.
\end{lem}

\begin{proof}
It will be enough to show that the canonical divisors 
$i^*(\Delta(\pi_* \lambda))$ and $\tau(D(\pi_* \lambda))$ are equal for a
general element $\lambda \in P_{odd}$. In
both cases the divisor coincide with the divisor $Nm(\delta)$,
where $\delta$ is the unique effective divisor in the linear system
$|\lambda|$. The computations are straight-forward and left to 
the reader. 
\end{proof}


\bigskip
\noindent
Therefore we can conclude that 
$$H^0(JC, \Lalpha)_{0-}^{(3)} = \ker \tau = \ker i^* =
\Gamma_{\tilde{C}}.$$
where $H^0(JC, \Lalpha)_{0-}^{(3)}$ denotes the
subspace of $H^0 ( JC,\Lalpha )_-$ of elements with multiplicity $\geq
2$ (hence $\geq 3$ by anti-symmetry) at the origin.  We now proceed as
in the proof of part 1 of Theorem 2. We consider the two linear maps
$$\phi_1 :  \Gamma_{\tilde{C}} \lra \Lambda^2 H^0(\omega \alpha)$$
$$\phi_2 : H^0(JC, \Lalpha)_{0-}^{(2)} \lra H^0(C \times C, \delta^*
\Lalpha (-2\Delta))_- = \Lambda^2 H^0(\omega \alpha)$$
which are defined as follows. As in part 1, $\phi_2$ is the map given by
pull-back under the difference map $\delta$. To define $\phi_1$, let
$N_{\tilde{C}/P'}$ denote the normal bundle of $i (\tilde{C})$ in
$P'$.  Then $\phi_1$ is obtained by restricting a section $s \in
\Gamma_{\tilde{C}}$ to the first infinitesimal neighborhood of
$\tilde{C}$.  In other words
$$\Gamma_{\tilde{C}}^{(2)} =  \ker \{ \phi_1 : \Gamma_{\tilde{C}} \lra
H^0(\tilde{C}, N_{\tilde{C}/P'}^* \otimes i^* \cO(2\Xi'_0))_- = 
H^0(\tilde{C}, N_{\tilde{C}/P'}^* \otimes \omega_{\tilde{C}})_- \} $$
The vector bundle $N_{\tilde{C}/P'}^*$ fits into the exact sequence
\begin{equation}\label{exseqN}
0\lra N_{\tilde{C}/P'}^*\lra H^0 (\omega\alpha )\otimes\cO_{\tC
  }\lra\omega_{\tC }\lra 0
\end{equation}
where the right-hand map is the embedding $H^0 (\omega\alpha )\otimes\cO_{\tC
  }\inj H^0 (\omega_{\tC } )\otimes\cO_{\tC
  }$ followed by evaluation $H^0 (\omega_{\tC } )\otimes\cO_{\tC
  }\ra\omega_{\tC }$. Therefore this map is the pull-back of
evaluation $H^0(\omega \alpha) \otimes \cO \stackrel{ev}{\ra}  \omega
\alpha$. Let $M$ be the kernel of the latter, i.e., we have
the exact sequence
\begin{equation}\label{exseqM}
0 \lra M \lra H^0(\omega \alpha) \otimes \cO \map{ev}  \omega \alpha 
\lra 0,
\end{equation}
whose pull-back by $\pi$ is \eqref{exseqN}.

We twist \eqref{exseqM} by $\omega \alpha$ and take cohomology
$$ 0 \lra H^0(C, M \otimes \omega \alpha) \lra H^0(\omega \alpha) \otimes
H^0(\omega \alpha) \map{m} H^0(\omega^2) \lra \ldots $$
where $m$ is the multiplication map. 
We deduce that 
$$H^0(\tilde{C}, N_{\tilde{C}/P'}^* \otimes \omega_{\tilde{C}})_- =
H^0(C,M \otimes \omega \alpha) = \ker m = \Lambda^2 H^0(\omega \alpha)
\oplus I_C^{Pr}(2)$$
where $I_C^{Pr}(2)$ is the space of quadrics through the Prym-canonical 
curve. It remains to show that $\im \phi_1 = \Lambda^2 H^0(\omega
\alpha)$. This will follow from the next two lemmas. First we will
compute, as in part 1, the image under $\phi_1$ of some special
sections $s_E \in \Gamma_{\tC}$, namely $s_E$ such that 
$Z(s_E) = \Delta(E)$ with
$E$ a general bundle in $\modkalpha$ with $h^0(E) = 2$, i.e.,
we determine the tangent spaces to $\Delta(E)$ along the
curve $i(\tC)$. This is done in the following lemma.

\begin{lem}\label{lemab}
Let $a,b$ be the sections defined by \eqref{defab}. Then we have 
$$\phi_1(s_E) = a \wedge b \in \Lambda^2 H^0(\omega \alpha)$$
up to multiplication by a nonzero scalar.
\end{lem}

\begin{proof}
First we need to show that for a general semi-stable
bundle $E$ with $h^0(E) =2$ the divisor $\Delta(E)$ is smooth
at a general point $i(\tp) \in \Delta(E)$. For this
decompose a general Prym-canonical divisor into two
effective divisors of degree $g-1$, i.e., $D + D' \in 
|\omega \alpha|$. Put $L = \cO(D)$. Then $h^0(D) = 1 =
h^0(\omega(-D)) = h^0(\omega \alpha(-D)) = h^0(\alpha(D))$.
If $E = L \oplus \omega \alpha L^{-1}$, then $\tE = \pi^*E =
\pi^*L \oplus \omega_{\tC} \pi^*L^{-1}$, $D(\tE) = 
\tT_{\pi^*L} + \tT_{\omega_{\tC}\pi^*L^{-1}}$ and
$\Delta(E) = 
\tT_{\pi^*L|P'} + \tT_{\omega_{\tC}\pi^*L^{-1}|P'}$.
At a general point $i(\tp) \in \tT_{\pi^*L}$, we
see immediately that the tangent space to $\tT_{\pi^*L}$ does
not contain the tangent space to $P'$, i.e., $\Delta(E)$ is 
smooth at $i(\tp)$. 
  Next we compute the tangent space to the divisor $\Delta(E)$ at a
  smooth point $i(\tp) \in \Delta(E)$.
The smoothness of $\Delta(E)$ at $i(\tp)$ implies that 
 $h^0(\tC, \tE(\tp - \sigma \tp)) = 2$. We choose a basis $\{
  u,v \}$ of the $2$-dimensional vector space $H^0(\tC,\tE(\tp -
  \sigma \tp))$. Then by \cite{laszlo} prop. V.2 and the same
  reasoning as in the proof of part 1 of Theorem 2, we see that the
  projectivized 
 tangent space $\TT_{i(\tp)}\Delta(E)$ to $\Delta(E)$ at $i(\tp)$,
  which is a hyperplane in $\pp T_{i(\tp)}P' \cong |\omega \alpha|^*$
  is the zero locus of the section in $\gamma(\tp) \in H^0(\omega \alpha)$,
 which is the  image of 
  $u \wedge \sigma^* v := u \otimes \sigma^*v - v \otimes \sigma^*u$
  under the exterior product map
  $$ H^0(\tE(\tp - \sigma \tp)) \otimes \sigma^*  
    H^0(\tE(\tp - \sigma \tp)) = H^0(\tE(\tp - \sigma \tp))
    \otimes H^0(\tE(\sigma \tp - \tp)) \map{\mu} H^0(\omega_{\tC}) $$
Since $\det E = \omega \alpha$, we see that $\gamma(\tp) =
\mu(u \wedge \sigma^* v) \in H^0(\omega \alpha) \subset 
H^0(\omega_{\tC})$. We will now describe the map
$\gamma: \tC \ra |\omega \alpha| \ : \ \tp \mapsto \gamma(\tp)$.
  Note that, since $h^0
  (\tE) =4$, we have $h^0 (\tE ( -\sigma\tp )) = 2$ for $\tp$ general.
  Hence $\{ u,v\}$ is also a basis for $H^0 (\tE (- \sigma\tp ))$. 
Consider the inclusion
$$H^0 (\tE (- \sigma\tp )) \subset H^0(\tE) = H^0(E) \oplus 
H^0(E\alpha)$$
and decompose $u = u_+ + u_-$, $v = v_+ + v_-$ with
$u_+,v_+ \in H^0(E) = H^0(\tE)_+$ and
$u_-,v_- \in H^0(E\alpha) = H^0(\tE)_-$. Then the element
$\gamma(\tp)$ is the image of
$(u_+ \wedge v_+, -u_- \wedge v_-) \in \Lambda^2 H^0(E) \oplus
\Lambda^2 H^0(E\alpha)$ under the exterior product map 
$\Lambda^2 H^0(E) \oplus \Lambda^2 H^0(E\alpha) \ra H^0(\omega 
\alpha)$, i.e., $\gamma(\tp) \in \pp(\cc a \oplus \cc b) \subset
|\omega \alpha|$. Since $\tC \subset \Delta(E)$, we have
$\varphi_{\alpha can}(p) \in \TT_{i(\tp)}(\Delta(E))$. So
for general $\tp$, $\gamma(\tp)$ is the unique
divisor of the pencil $\pp(\cc a \oplus \cc b)$ containing
$\tp$. Hence we can conclude  
 that the section $\phi_1(s_E) \in H^0(M \otimes \omega \alpha)$
considered as a tensor in $H^0(\omega \alpha) \otimes H^0(\omega
\alpha)$ is $a \wedge b$.
\end{proof}

Since, a priori, we do not know that $\pp \Gamma_{\tC}$ is spanned
by divisors of the form $\Delta(E)$, we need to establish
a symmetry property for any divisor $D \in \pp \Gamma_{\tC}$. This
is done as follows.

Let $\ts,\tt \in \tC$ be two points of $\tC$ with respective images
$s,t \in C$ and let $D$ be an element of $\pp \Gamma_{\tC}$.  Assume
that $i(\ts),i(\tt) \in D$ are smooth points of $D$ and let $\TT_s D$
and $\TT_t D$ denote the projectivized tangent spaces to the divisor $D$
at the points $i(\ts)$ and $i(\tt)$. Since we can identify the
projectivized tangent space to the Prym variety $P'$ at any point with
the Prym-canonical space $|\omega \alpha|^*$, we may view $\TT_s D$
and $\TT_t D$ as hyperplanes in $|\omega \alpha|^*$.  Note that $\TT_s
D$ only depends on $s \in C$ and not on the lift $\ts \in \tC$. Then
we have

\begin{lem}\label{lemst}
With the preceding notation, we have an equivalence
$$ \varphi_{\alpha can}(s) \in \TT_t D \ \iff \ 
 \varphi_{\alpha can}(t) \in \TT_s D $$
\end{lem}

\begin{proof}
  Consider the invertible sheaf $x = \cO_{\tC}(\ts - \sigma \ts + \tt
  - \sigma \tt) \in P$ and the corresponding embedding
$$
i_x : \tC \lra P' \qquad \tp \lms \cO_{\tC}(\tp - \sigma \tp)
\otimes x. $$
The curve $i_x(\tC)$ is the curve $i(\tC)$ translated by
$x$. A straight-forward computation shows that $i_x^{ -1 }(\cO_{P'}
(2\Xi'_0)) = \omega_{\tC}x^{-2}$ and by a result of Beauville (see
\cite{izadivanstraten} page 569) the induced linear map on global
sections $H^0(P',2\Xi'_0) \ra H^0(\omega_{\tC}x^{-2})$ is surjective.
We observe that
\[
i_x(\sigma \tt) = i(\ts), \qquad i_x(\sigma \ts) = i(\tt),
\]
and that the projectivized tangent line to the curve $i_x(\tC)$ at the
point $i_x(\sigma \tt)$ (resp. $i_x(\sigma \ts)$) is the point
$\varphi_{\alpha can}(t)$ (resp. $\varphi_{\alpha can}(s)$) in
$|\omega \alpha|^* \cong \pp T_{i(\ts)}P'$ (resp. $\cong \pp
T_{i(\tt)}P'$).  Let $\TT_{\ts}$ (resp. $\TT_{\tt}$) denote the
embedded tangent line in $| 2\Xi_0'|^*$ to the curve $i_x(\tC)$ at the
point $i_x(\sigma \tt)$ (resp. $i_x(\sigma \ts)$), so that $\TT_{\ts}$
(resp.  $\TT_{\tt}$) passes through the point $i(\ts)$ (resp.
$i(\tt)$) with tangent direction $\varphi_{\alpha can}(t)$ (resp.
$\varphi_{\alpha can}(s)$). Then the lemma will follow if we show that
these two tangent lines intersect in a point $I(\ts,\tt)$, i.e.
\begin{equation} \label{intpoint}
\TT_{\ts} \cap \TT_{\tt} = I(\ts,\tt) \in |2\Xi'_0|^*.
\end{equation}
This property follows from a dimension count: since $C$ is
non-hyperelliptic, we have $x^{-2} \not= \cO_{\tC}$, so
$h^0(\omega_{\tC}x^{-2}) = 2g-2$. Since $h^0(\omega_{\tC}x^{-2}
(-2\sigma \ts - 2\sigma \tt)) = h^0(\omega_{\tC}(-2\ts - 2\tt))
\geq 2g-5$, the tangent lines $\TT_{\tt}$ and
$\TT_{\ts}$ are contained in a projective $2$-plane,
hence intersect. To get the equivalence stated in the lemma, let $H_D$
denote the hyperplane in $|2\Xi'_0|^*$ corresponding to the
divisor $D \in \pp \Gamma_{\tC}$. Assume e.g. that 
$\varphi_{\alpha can}(s) \in \TT_t D$. This means that $H_D$ contains
$\TT_{\tt}$. Since $i(\ts) \in H_D$, it follows from
\eqref{intpoint} that $H_D$ also contains 
$\TT_{\ts}$, so $\varphi_{\alpha can}(t) \in \TT_s D$. 
\end{proof}

At this stage we can conclude: by lemma \ref{lemst} we know that for
all $s \in  \Gamma_{\tC}$, $\phi_1(s) \in H^0(\omega \alpha)
\otimes H^0(\omega \alpha)$ lies either in the symmetric or
skew-symmetric eigenspace, i.e. $\im \phi_1 \subset I_C^{Pr}(2)
\subset \sym^2 H^0(\omega \alpha)$ or $\im \phi_1 \subset \Lambda^2
H^0(\omega \alpha)$. Lemma \ref{lemab} asserts that $\im \phi_1
\subset \Lambda^2 H^0(\omega \alpha)$.

As in \eqref{comprdelta}, we have that $\phi_2(pr_-(s_E)) =
pr_-(\delta^*(s_E))$, where $pr_-$ denotes the projection $ H^0(\omega
\alpha) \otimes H^0(\omega \alpha) \lra \Lambda^2 H^0(\omega \alpha)$
and $s_E$ is as above. Hence we see that $\phi_2(pr_-(s_E)) = a \wedge
b$. By lemma \ref{lemab} the projectivizations of $\phi_1$ and
$\phi_2$ coincide on all divisors of the form $\Delta(E)$ whose images
generate $\pp \wedge^2 H^0 (\omega\alpha)$. Hence $\phi_1 = \phi_2$ up to
a nonzero scalar and $\phi_1$ and $\phi_2$ are surjective.

\begin{rem}
  An alternative way of proving that $\im \phi_1 \subset \Lambda ^2
  H^0(\omega \alpha)$ would be to twice take the derivative of the
  quadrisecant identity for Prym varieties \cite{fay} prop. 6 (fix two
  points and consider the other two as canonical coordinates on the
  universal cover of $\tC$.)
\end{rem}

\section{The scheme-theoretical base locus of $\pp \Gamma_{\tC}$}

>From section 2 we know that the sets $Bs(\pp \Gamma_{\tC})$ and $i(\tC)$
are equal. To prove the scheme-theoretical equality, it will be
enough to show that, $\forall \tp \in \tC$, the projectivized
tangent spaces at $i(\tp)$ to divisors $D \in \pp \Gamma_{\tC}$ cut 
out the projectivized tangent space at $i(\tp)$ to $i(\tC)$, 
which is $\varphi_{\alpha can}(p) \in |\omega \alpha|^* = \pp
T_{i(\tp)}P'$, i.e.,
\begin{equation} \label{inttan}
\bigcap_{D \in \pp \Gamma_{\tC}} \TT_{i(\tp)} D = \varphi_{\alpha can}(p).
\end{equation}
If we take $D = \Delta(E)$ for some semi-stable vector bundle
$E$ with $h^0(E) = 2$ (see section 4.2) then the hyperplane
$\TT_{i(\tp)}(\Delta(E)) \subset |\omega \alpha|^*$
corresponds to the unique section of the pencil $\pp(\cc a \oplus
\cc b)$ vainshing at $p$ (proof of lemma 4.3). Since for general
$a,b \in |\omega \alpha|$ we can find a vector bundle $E$ (lemma 4.1)
such that equality in lemma 4.3 holds, we can conclude \eqref{inttan}.

\flushleft{Elham Izadi \\
 Department of Mathematics \\ 
 Boyd Graduate Studies Research Center \\
 University of Georgia \\ 
 Athens, GA 30602-7403 \\
 USA \\
 email: izadi@math.uga.edu}
\bigskip
\noindent
\flushleft{Christian Pauly \\
Laboratoire J.-A. Dieudonn\'e \\
Universit\'e de Nice Sophia Antipolis \\
Parc Valrose \\
06108 Nice Cedex 02 \\ France \\
email: pauly@math.unice.fr}

\end{document}